\documentclass[11pt]{amsart}
\usepackage[hyperindex]{hyperref} 
\usepackage{amsrefs}
\title{The Laguerre polynomials preserve real-rootedness}
\author{Steve fisk\\Bowdoin College\\Brunswick Maine 04011}

\newtheorem{theorem}{Theorem}

\begin{document}
\maketitle

The study of linear transformations that map polynomials with all real
roots to polynomials  with all real roots has been of interest for 
many years \cites{me,polya-szego}. We say such transformations \emph{preserve
  real-rootedness}. Notice that a linear transformation
acting on polynomials can be defined by specifying its action on the polynomials
$x^n$.

\begin{theorem}
  If $L_n(x)$ is the $n$'th Laguerre polynomial then the linear
  transformation $x^n\mapsto L_n(x)$ preserves real-rootedness.
\end{theorem}
\begin{proof}
  We recall a basic result from \cite{polya-szego}*{Part V, No. 65}: if $T(x^k)=x^k/k!$  then $T$
  preserves real-rootedness. The Laguerre polynomials satisfy
  \cite{szego}*{page 101}
\[
L_n(x) = \sum_{i=0}^n (-1)^k\binom{n}{k}\frac{x^k}{k!}  
\] 
which equals $T(1-x)^n$. Thus the linear transformation $x^n\mapsto
L_n(x)$ is the composition of two linear transformations preserving
real-rootedness, namely $x\mapsto 1-x$ and $x^k\mapsto x^k/k!$, and so
preserves real-rootedness.
\end{proof}

\end{document}